\documentclass[12pt]{amsart}
\usepackage{amssymb,amsmath}
\def\nova{u^{\tau + (0)}}
\def\supernova{u^{\tau + (H)}}
\def\Op{\T{Op}}
\def\Tan{\T{Tan}}
\def\NO#1{||#1||^2}
\def\no#1{||#1||}
\def\ni#1{|||#1|||}
\def\nno#1{|||#1|||}
\numberwithin{equation}{section}
\def\db{\bar\partial}
\def\db*{\bar\partial^*}
\def\T{\text}
\def\simgeq{\underset\sim>}
\def\simleq{\underset\sim<}
\def\1#1{\overline{#1}}
\def\2#1{\widetilde{#1}}
\def\3#1{\widehat{#1}}
\def\4#1{\mathbb{#1}}

\def\5#1{\frak{#1}}
\def\6#1{{\mathcal{#1}}}
\def\C{{\4C}}
\def\R{{\4R}}

\def\Z{{\4Z}}

\def\di{\partial}
\def\dib{\bar\partial}
\begin{document}
\title[regularity at the boundary and tangential regularity...]{regularity at the boundary and tangential regularity}
\author[T.V.~Khanh and G.~Zampieri]{Tran Vu Khanh and  Giuseppe Zampieri}
\maketitle
\begin{abstract}
For a pseudoconvex domain $D\subset \C^n$, we prove the equivalence of the local hypoellipticity of the system $(\dib,\dib^*)$ with the system $(\dib_b,\dib^*_b)$ induced in the boundary. This develops our former result in \cite{KZ10} which used the theory of the ``harmonic" extension by Kohn. This technique is inadequate for the purpose of the present paper and must be replaced by the ``holomorphic" extension introduced by the authors in \cite{KZ09}.
\end{abstract}
\def\Giialpha{\mathcal G^{i,i\alpha}}
\def\cn{{\C^n}}
\def\cnn{{\C^{n'}}}
\def\ocn{\2{\C^n}}
\def\ocnn{\2{\C^{n'}}}
\def\const{{\rm const}}
\def\rk{{\rm rank\,}}
\def\id{{\sf id}}
\def\aut{{\sf aut}}
\def\Aut{{\sf Aut}}
\def\CR{{\rm CR}}
\def\GL{{\sf GL}}
\def\Re{{\sf Re}\,}
\def\Im{{\sf Im}\,}
\def\codim{{\rm codim}}
\def\crd{\dim_{{\rm CR}}}
\def\crc{{\rm codim_{CR}}}
\def\phi{\varphi}
\def\eps{\varepsilon}
\def\d{\partial}
\def\a{\alpha}
\def\b{\beta}
\def\g{\gamma}
\def\G{\Gamma}
\def\D{\Delta}
\def\Om{\Omega}
\def\k{\kappa}
\def\l{\lambda}
\def\L{\Lambda}
\def\z{{\bar z}}
\def\w{{\bar w}}
\def\Z{{\1Z}}
\def\t{{\tau}}
\def\th{\theta}
\emergencystretch15pt
\frenchspacing
\newtheorem{Thm}{Theorem}[section]
\newtheorem{Cor}[Thm]{Corollary}
\newtheorem{Pro}[Thm]{Proposition}
\newtheorem{Lem}[Thm]{Lemma}
\theoremstyle{definition}\newtheorem{Def}[Thm]{Definition}
\theoremstyle{remark}
\newtheorem{Rem}[Thm]{Remark}
\newtheorem{Exa}[Thm]{Example}
\newtheorem{Exs}[Thm]{Examples}
\def\Label#1{\label{#1}}
\def\bl{\begin{Lem}}
\def\el{\end{Lem}}
\def\bp{\begin{Pro}}
\def\ep{\end{Pro}}
\def\bt{\begin{Thm}}
\def\et{\end{Thm}}
\def\bc{\begin{Cor}}
\def\ec{\end{Cor}}
\def\bd{\begin{Def}}
\def\ed{\end{Def}}
\def\br{\begin{Rem}}
\def\er{\end{Rem}}
\def\be{\begin{Exa}}
\def\ee{\end{Exa}}
\def\bpf{\begin{proof}}
\def\epf{\end{proof}}
\def\ben{\begin{enumerate}}
\def\een{\end{enumerate}}
\def\dotgamma{\Gamma}
\def\dothatgamma{ {\hat\Gamma}}

\def\simto{\overset\sim\to\to}
\def\1alpha{[\frac1\alpha]}
\def\T{\text}
\def\R{{\Bbb R}}
\def\I{{\Bbb I}}
\def\C{{\Bbb C}}
\def\Z{{\Bbb Z}}
\def\Fialpha{{\mathcal F^{i,\alpha}}}
\def\Fiialpha{{\mathcal F^{i,i\alpha}}}
\def\Figamma{{\mathcal F^{i,\gamma}}}
\def\Real{\Re}
%
%
%
\section{ }
Let $D$ be a pseudoconvex domain of $\C^n$  defined by $r<0$ with $C^\infty$ boundary $bD$. We use the standard notations $\Box=\dib\dib^*+\dib^*\dib$ for the complex Laplacian and $Q(u,u)=\no{\dib u}^2+\no{\dib^*u}^2$ for the energy form and some variants as, for an operator $\Op$, $Q_{\Op}(u,u)=\no{\Op \dib u}^2+\no{\Op\dib^*u}^2$. Here $u$ is a antiholomorphic form of degree $k\leq n-1$ belonging to $D_{\dib^*}$. We similarly define the tangential version of these objects, that is, $\Box_b,\,\,\dib_b,\,\,\dib^*_b,\,\,Q^b_{\Op}$. We take local coordinates $(x,r)$ in $\C^n$ with  $x\in\R^{2n-1}$ being the tangential coordinates and $r$, the equation of $bD$, serving as the last coordinate. We define the tangential $s$-Sobolev norm by $\ni{u}_s:=\no{\Lambda^su}_0$ where $\Lambda^s$ is the standard tangential pseudodifferential operator with symbol $\Lambda^s_\xi=(1+|\xi|^2)^{\frac s2}$.
We note that
\begin{equation}
\Label{1}
\begin{cases}
\NO{\dib u}_s+\NO{\dib^*u}_s=\underset{j\leq s}\sum Q_{\Lambda^{s-j}\di_r^j}(u,u),
\\
|||\dib u|||^2_s+|||\dib^* u|||^2_s=Q_{\Lambda^s}(u,u),
\\
\NO{\dib_b u_b}_s+\NO{\dib^*_bu_b}_s=Q^b_{\Lambda^s}(u_b,u_b).
\end{cases}
\end{equation}
We decompose $u$ in tangential and normal component, that is
$$
u=u^\tau+u^\nu,
$$
and further decompose in microlocal components (cf. \cite{K02})
$$
u^\tau=u^{\tau+}+u^{\tau-}+u^{\tau0}.
$$
We similarly decompose $u_b=u_b^++u_b^-+u_b^0$. 
We use the notation $\bar L_n$ for the ``normal" $(0,1)$-vector field and $\bar L_1,...,\bar L_{n-1}$ for the tangential ones. We have therefore the description 
for the totally real tangential, resp. normal, vector field $T$, resp $\di_r$:
\begin{equation*}
\begin{cases}
T=i(L_n-\bar L_n),
\\
\di_r=L_n+\bar L_n.
\end{cases}
\end{equation*}
From this, we get back $\bar L_n=\frac12(\di_r+iT)$. We denote by $\sigma$ the symbol of a (pseudo)differential operator and  by $\tilde u$ the partial tangential Fourier transform of $u$. We define a  ``holomorphic" extension 
$u^{\tau+(H)}$ by 
\begin{equation}
\Label{1.-1}
\supernova=(2\pi)^{-2n+1}\int_{\R^{2n-1}}e^{ix\xi}e^{r\sigma(T)}\psi^+(\xi)\tilde u(\xi,0)d\xi.
\end{equation}
This definition has been introduced in \cite{KZ09}.
Note that  $\sigma(T)\simgeq (1+|\xi|^2)^{\frac12}$ for $\xi$ in supp$\,\psi^+$ and $(x,r)$ in a local patch; thus in the integral, the exponential is dominated by $e^{-|r|(1+|\xi|^2)^{\frac12}}$ for $r<0$. 
Differently from the harmonic extension by Kohn, the present one is well defined only in positive microlocalization.
 We can think of  $\supernova$ in two different ways: either as a modification of $u^{\tau+}$ or as an extension of $u_b^+$.
We have a first relation from \cite{K02} p. 241, between a trace $v_b$ and a general extension $v$: for any $\epsilon$ and suitable $c_\epsilon$ 
\begin{equation}
\Label{1.5}
\no{v_b}\simleq c_\epsilon\ni{v}_{\frac12}+\epsilon\ni{\di_rv}_{\frac12}.
\end{equation}
This can been seen in \cite{K02} p. 241 and \cite{KZ09} as for the small/large constant argument.
As a specific property of our extension we have  the reciprocal relation to \eqref{1.5}, that is
\begin{equation}
\Label{1.7}
\no{r^k\supernova}\simleq \no{u^+_b}_{-k-\frac12}.
\end{equation}
This is readily checked ( \cite{KZ09} (1.12)). 
We denote by $\bar\partial^\tau$ the extension of $\bar\partial_b$ from $b\Omega$ to $\Omega$ which stays tangential to the level surfaces $r\equiv\T{const}$. It acts on tangential forms $u^\tau$ and its action  is  $\bar\partial^\tau u^\tau=(\bar\partial u^\tau)^\tau$. We denote by $\bar\partial^{\tau\,*}$ its adjoint; thus $\bar\partial^{\tau\,*}u^\tau=\bar\partial^*(u^\tau)$. We use the notations $\Box^\tau$ and $Q^\tau$ for the corresponding Laplacian and energy form.
 We notice that 
\begin{equation} 
\Label{2.3ter}
\begin{split}
Q(\supernova,\supernova)&=Q^\tau(\supernova,\supernova)+\NO{\bar L_n \supernova}_0
\\
&=Q^\tau(\supernova,\supernova).
\end{split}
\end{equation}
We have to describe how \eqref{1.5} and \eqref{1.7} are affected by $\dib$ and $\dib^*$.
\bp
\Label{p1.1}
We have for any extension $v$ of $v_b$
\begin{equation}
\Label{1.3}
Q^b(v_b,v_b)\simleq Q^\tau_{\Lambda^{\frac12}}(v,v)+Q^\tau_{\di_r\Lambda^{-\frac12}}(v,v),
\end{equation}
and, specifically for $\supernova$
\begin{equation}
\Label{1.4}
\begin{cases}
Q^\tau(\supernova,\supernova)\simleq Q^b_{\Lambda^{\frac12}}(u^+_b,u^+_b)+\NO{u^+_b}_{-\frac12}
\\
\bar L_n\supernova\equiv0.
\end{cases}
\end{equation}
\ep
\bpf
We have
$$
\dib^\tau v|_{bD}=\dib_bv_b,\qquad\dib^{\tau*}v|_{bD}=\dib^*_bv_b.
$$
Then, \eqref{1.3} follows from \eqref{1.5}.

We pass to prove \eqref{1.4}. We have $\dib^\tau=\dib_b+r\Tan$, $\dib^{\tau*}=\dib^*_b+r\Tan$ which yields
\begin{equation}
\Label{1.6}
\begin{cases}
\dib^\tau \supernova=(\dib_bu_b)^{\tau+(H)}+r\Tan\, \supernova,
\\
\dib^{\tau*}\supernova=(\dib^*_bu_b)^{\tau+(H)}+r\Tan\,\supernova.
\end{cases}
\end{equation}
Application of \eqref{1.7} yields
\begin{equation*}
\begin{split}
\NO{\dib^\tau\supernova}+\NO{\dib^{\tau*}\supernova}&=\NO{(\dib_bu_b)^{\tau+(H)}}+\NO{(\dib^*_bu_b)^{\tau+(H)}}+\NO{r\Tan \,\supernova}
\\
&\simleq \NO{\dib_bu^+_b}_{-\frac12}+\NO{\dib^*_bu^+_b}_{-\frac12}+\NO{u^+_b}_{-\frac12},
\end{split}
\end{equation*}
which is the first of \eqref{1.4}. The second is an easy consequence of the relation $\bar L_n=\frac12(\di_r+iT).$ 

\epf
We finally decompose $ u^{\tau+}=\supernova+\nova$ which also serves as a definition of $\nova$.

\bp
\Label{p1.2}
Each of the forms $u^\#=u^\nu,\,\,u^{\tau\,-},\,\,u^{\tau\,0},\,\,\nova,\,\,u^-_b,\,\,u^0_b$ enjoys elliptic estimates, that is 
\begin{equation}
\Label{1.8}
\no{\zeta u^\#}_s\simleq\no{\zeta'\bar\partial u^\#}_{s-1}+\no{\zeta'\bar\partial^* u^\#}_{s-1}+\no{u^\#}_0\qquad s\geq2.
\end{equation}
\ep
\bpf
Estimate \eqref{1.8} follows, by iteration, from
\begin{equation}
\Label{1.8bis}
\no{\zeta u^\#}_s\simleq\no{\zeta\bar\partial u^\#}_{s-1}+\no{\zeta\bar\partial^* u^\#}_{s-1}+\no{\zeta'u^\#}_{s-1}.
\end{equation}
As for $u^\nu$ and $\nova$ this latter follows from $u^\nu|_{bD}\equiv0$ and $\nova|_{bD}\equiv0$. For the terms with $-$ and $0$, this
 follows from the fact that $|\xi_T|\simleq |\sigma(\dib)|$ in the region of $0$-micolocalization and from $\sigma[\dib,\dib^*]\leq0$ and $\sigma(T)<0$ in the negative microlocalization. 
 We refer to \cite{FK72} formula (1) of Main theorem as a general reference but also give an outline of the proof. 
We start from
\begin{equation}
\Label{2.3bis}
\nno{\zeta u^\#}^2_1\simleq Q(\zeta u^\#,\zeta u^\#)+\no{\zeta'u^\#}^2_0;
\end{equation}
this is the basic estimate for $u^\nu$ and $\nova$ (which vanish at $bD$) whereas it is \cite{K02} Lemma 8.6 for $u^{\tau\,-}$, $u^{\tau\,\,0}$ and $u^-_b$, $u^0_b$. Applying 
\eqref{2.3bis} to $\zeta\Lambda^{s-1}\zeta u^\#$ one gets the estimate of tangential norms for any $s$, that is, \eqref{1.8bis} with the usual norm replaced by the ``triplet" norm.  Finally, by non-characteristicity of $(\bar\partial,\bar\partial^*)$ one passes from tangential to full norms along the guidelines of \cite{Z08} Theorem~1.9.7. The version of this argument for $\Box$ can be found in \cite{K02} second part of p. 245.

\epf
Let $\zeta$ and $\zeta'$ be a couple of cut-off with $\zeta\prec\zeta'$ in the sense that $\zeta'|_{\T{supp}\,\zeta}\equiv1$, and let $s$ and $l$ be a pair of indices.
\bt
\Label{t1.1}
The following two estimates are equivalent:
\begin{gather}
\Label{1.9}
\no{\zeta u_b}_s\simleq\no{\zeta'\dib_bu_b}_{s+l}+\no{\zeta'\dib^*_bu_b}_{s+l}+\no{u_b}_0\quad \T{ for any $u_b\in C^\infty_c(b\Om\cap U)$},
\\
\Label{1.10}
\no{\zeta u}_s\simleq\no{\zeta'\dib u}_{s+l}+\no{\zeta'\dib^*u}_{s+l}+\no{u}_0\quad\T{for any $u\in D_{\dib^*}\cap C^\infty_c(\bar \Om\cap U)$}.
\end{gather}
\et
\br
The above estimates \eqref{1.9} and \eqref{1.10} for any $s,\,\zeta,\,\,\zeta'$ and for suitable $l$, characterize the local hypoellipticity of the system $(\dib_b,\dib^*_b)$ and $(\dib,\dib^*)$ respectively (cf. \cite{K05}). When $l>0$, one says that the system has a ``loss" of $l$ derivatives; when $l<0$, one says that it has a ``gain" of $-l$ derivatives. 
\er
\bpf
Because of Proposition~\ref{p1.2}, it suffices to prove \eqref{1.9} for $u^+_b$ and \eqref{1.10} for $u^{\tau\,+}$. It is also obviuos that we can consider cut-off functions $\zeta$ and $\zeta'$ in the only tangential coordinates, not in $r$.
 We start by proving that \eqref{1.9} implies \eqref{1.10}. We recall the decomposition $u^{\tau+}=\supernova+u^{\tau+(0)}$ and begin by estimating $\supernova$. We have
\begin{equation}
\Label{1.12}
\begin{split}
\ni{\zeta\supernova}^2_s&\underset{\eqref{1.7}}\simleq \NO{\zeta u_b^+}_{s-\frac12}
\\
&\underset{\eqref{1.9}}\simleq Q^b_{\Lambda^{s+l-\frac12}\zeta'}(u^+_b,u^+_b)+\NO{u^+_b}_{-\frac12}
\\
&\underset{\eqref{1.3}}\simleq Q^\tau_{\Lambda^{s+l}\zeta'}(u^{\tau+},u^{\tau+})+Q^\tau_{\di_r\Lambda^{s+l-\frac12}\zeta'}(u^{\tau+},u^{\tau+})+\NO{u^{\tau+}}_0.
\end{split}
\end{equation}
It remains to estimate $u^{\tau+(0)}$; 
since $u^{\tau+(0)}|_{bD}\equiv0$, then by $1$-elliptic estimates
\begin{equation}
\Label{2}
\begin{split}
\ni{\zeta u^{\tau+(0)}}_s&\underset{\eqref{1.8bis}}\simleq Q_{\Lambda^{s-1}\zeta}( \nova,\nova)+\ni{\zeta'\nova}_{s-1}^2
\\
&\simleq Q_{\Lambda^{s-1}\zeta}(u^{\tau+},u^{\tau+})+Q^\tau_{\Lambda^{s-1}\zeta}(\supernova,\supernova)+\ni{\zeta'\nova}_{s-1}^2
\\
&\simleq Q_{\Lambda^{s-1}\zeta}(u^{\tau+},u^{\tau+})+\ni{\zeta \supernova}_s^2+\ni{\zeta'u^{\tau+(H)}}_{s-1}^2+\ni{\zeta'\nova}_{s-1}^2,
\end{split}
\end{equation}
where we have used that $Q=Q^\tau$ over $\supernova$ in the second inequality together with the estimate $Q^\tau_{\Lambda^{s-1}}\simleq \Lambda^s$ in the third. 
We estimate terms in the last line.
First, the term $\ni{\zeta \supernova}_{s}^2$ 
is estimated by means of \eqref{1.12}. Next, the terms in $(s-1)$-norm can be brought to $0$-norm by combined inductive use of \eqref{1.12} and \eqref{2} and eventually their sum is controlled by $\no{ u^{\tau+}}_0^2$. We put together \eqref{1.12} and \eqref{2} (with the above further reductions), recall the first of \eqref{1} in order to estimate $Q^\tau{\Lambda^s\zeta'}+Q^\tau_{\di_r\Lambda^{s-1}\zeta'}$ in the right of \eqref{1.12} and end up with
\begin{equation}
\Label{3}
\ni{\zeta u^{\tau+}}_s\simleq \no{\zeta'\dib u^{\tau+}}_s+\no{\zeta'\dib^* u^{\tau+}}_s+\no{u^{\tau+}}_0.
\end{equation}
Finally, by non-characteristicity of $(\bar\partial,\bar\partial^*)$ one passes from tangential to full norms in the left side of \eqref{3} along the guidelines of \cite{Z08} Theorem~1.9.7. The version of this argument for $\Box$ can be found in \cite{K02} second part of p. 245. Thus we get \eqref{1.10}.

We prove the converse. Thanks to $\di_r=\bar L_n+\Tan$ and to the second of \eqref{1.4}, we have $\di_r\supernova=\Tan\,\supernova$.
It follows
\begin{equation}
\Label{1.13}
\begin{split}
\NO{\zeta u^+_b}_s&\underset{\eqref{1.5}}\simleq c_\epsilon\ni{\zeta\supernova}^2_{s+\frac12}+\epsilon\ni{\di_r\zeta\supernova}^2_{s-\frac12}
\\
&\underset{\eqref{1.10}}\simleq Q^\tau_{\Lambda^{s+l+\frac12}\zeta'}(\supernova,\supernova)+\epsilon\ni{\zeta\supernova}_{s+\frac12}^2
\\
&\underset{\eqref{1.4}}\simleq Q^b_{\Lambda^{s+l}\zeta'}(u^+_b,u^+_b)+\epsilon\NO{\zeta u^+_b}_{s}.
\end{split}
\end{equation}
We absorb the term with $\epsilon$ and get  \eqref{1.9}.

\epf
Let $N$ and $G$ be the Neumann and Green operators, that is, the $H^0$-inverse of $\Box$ in $D$ and $\Box_b$ in $bD$ respectively.  
\br
\eqref{1.9} and \eqref{1.10} imply local regularity, but not exact $s$-Sobolev regularity, of $G$ and $N$ respectively. We first prove for $N$. We start from remarking that
\begin{equation}
\Label{6}
\begin{cases}
\dib^*N\T{ is exactly regular over $\T{Ker} \dib$},
\\
\dib N\T{ is exactly regular over $\T{Ker} \dib^*$}.
\end{cases}
\end{equation}
As for the first, we put $u=\bar\partial^*Nf$ for $f\in\T{Ker}\,\bar\partial$. We have $(\dib u=f,\,\dib^*u=0)$ and hence by \eqref{1.10} $\no{\zeta u}_s\simleq \no{\zeta'f}_s+\no{u}_0$. To prove the second, we have just to put $u=\dib Nf$ for $f\in\T{Ker}\,\dib^*$ and reason likewise.
It follows from \eqref{6}, that  the Bergman projection $B$ is also  regular. 
(Notice that exact regularity is perhaps lost by taking the additional $\dib$ in $B:=\T{Id}-\dib^*N\dib$.) 
Finally, we exploit formula (5.36) in \cite{S10} in unweighted norms, that is, for $t=0$:
\begin{equation}
\Label{BS}
\begin{split}
N_q&=B_q(N_q\bar\partial)(Id-B_{q-1})(\db*N_q)B_q
\\&\quad +(Id-B_q)(\db*N_{q+1})B_{q+1}(N_{q+1}\bar\partial )(Id-B_q).
\end{split}
\end{equation}
Now, in the right side, the $\bar\partial  N$'s and $\db* N$'s are evaluated over $\T{Ker} \db*$ and $\T{Ker}\,\bar\partial $ respectively; thus they are exactly regular. The $B$'s are also regular and therefore such is $N$.
This concludes the proof of the regularity of $N$. The proof of the regularity of $G$ is similar, apart from replacing \eqref{BS} by its version for the Green operator $G$ stated in Section 5 of \cite{Kh10}.
\er

\hskip11cm $\Box$

\end{document}